\documentclass{amsart}

\usepackage{amsmath, amsthm, amssymb}
\newtheorem{thm}{Theorem}[section]
\newtheorem{cor}[thm]{Corollary}
\newtheorem{lem}[thm]{Lemma}

\newtheorem{rmk}[thm]{Remark}
\newtheorem{lemma}[thm]{Lemma}
\newtheorem{definition}[thm]{Definition}

\newcommand{\h}{{\bf h}}

\newcommand{\af}{\text{af}}
\newcommand{\fl}{\widetilde{Fl}}
\DeclareMathOperator{\Strong}{Strong}
\DeclareMathOperator{\Weak}{Weak}
\DeclareMathOperator{\outside}{outside}
\DeclareMathOperator{\inside}{inside}

\title{Pieri rule for the affine flag variety}
\author{Seung Jin Lee}
\address{School of Mathematics, Korea Institute for Advanced Study, 85 Hoegiro Dongdaemun-gu, Seoul 130-722, Republic of Korea.}
\email{lsjin@kias.re.kr}
\date{}

\begin{document}
\maketitle
\begin{abstract}
We prove the affine Pieri rule for the cohomology of the affine flag variety conjectured by Lam, Lapointe, Morse and Shimozono.
We study the cap operator on the affine nilHecke ring that is motivated by Kostant and Kumar's work on the equivariant cohomology of the affine flag variety. We show that the cap operators for Pieri elements are the same as Pieri operators defined by Berg, Saliola and Serrano. This establishes the affine Pieri rule.
\end{abstract}

\section{Introduction}

Affine Schubert calculus is a subject that ties combinatorics, algebraic geometry and representation theory together. Its modern development is motivated by the relation between $k$-Schur functions and the (co)homology of the affine Grassmannian of $SL(n)$. $k$-Schur functions were introduced by Lapointe, Lascoux, Morse \cite{LLM03} in the study of Macdonald polynomial positivity, a mostly combinatorial branch of symmetric function theory.
\\

 Peterson \cite{Pet97} identified the equivariant homology of the affine Grassmannian with a subalgebra of the affine nilHecke algebra $\mathbb{A}$, now called the Peterson algebra. Lam \cite{Lam08} proved that $k$-Schur functions can be identified with the Schubert classes of the homology of the affine Grassmannian of $SL(n)$. The nilHecke ring acts as divided difference operators on the equivariant cohomology of Kac-Moody partial flag varieties. By using the correspondence, one can investigate problems about the (co)homology of the affine Grassmannian of $SL(n)$ by translating them into the theory of $k$-Schur functions and understanding the combinatorics of $k$-Schur functions. There are affine analogues of the classical theory of Pieri rules, tableaux, and Stanley symmetric functions \cite{Lam06,LM05,LM07}. \\

Lam, Lapointe, Morse and Shimozono \cite{LLMS10} introduced the strong Schur functions, indexed by elements 
in the affine symmetric group $W_\af$. These strong Schur functions generalize the $k$-Schur functions combinatorially. They conjectured a series of properties of strong Schur functions such as the symmetry of the strong Schur functions. Berg, Saliola and Serrano \cite{BSS11,BSS12} studied the Pieri operators acting on the affine nilCoxeter ring $\mathbb{A}_0$ to establish some of conjectures in \cite{LLMS10}. \\

In this paper, we prove the affine Pieri rule for the cohomology of the affine flag variety conjectured in \cite{LLMS10}. We introduce the cap operators acting on the affine nilCoxeter ring $\mathbb{A}_0$ by investigating the work of Kostant and Kumar \cite{KK86} and show that the cap operators for Pieri elements are the same as the Pieri operators defined in \cite{BSS12} by using the strong strips. The affine Pieri rule gives us geometric interpretation of the skew strong Schur functions as an affine Grassmannian part of the cap product of the Schubert classes in (co)homology of the affine flag variety. We now describe these two operators.

\subsection{Pieri operators} 
For $u,v\in W_\af$, a \emph{marked strong cover} $C=\left(u\overset{a}\longrightarrow v \right)$ consists of $u,v\in W_\af$ and an integer $a$ such that $u=vt_{ij},v\lessdot u, a= u(j)=v(i)$ where $i\leq 0 <j$ and $\lessdot$ is the Bruhat order. We use the notation inside$(C)=u$ and outside$(C)=v$.\\
Lam, Lapointe, Morse and Shimozono \cite{LLMS10} introduced the notion of \emph{strong strips} to define $k$-Schur functions and strong Schur functions. A \emph{strong strip} $S$ of length $i$ from $u$ to $v$, denoted by $u\longrightarrow_{(i)} v$, is a path
\[u\overset{a_1}\longrightarrow u_1 \overset{a_2}\longrightarrow\cdots \overset{a_i}\longrightarrow u_i=v\]
where $a_1>a_2>\ldots>a_i$. We let inside$(S)=u$, outside$(S)=v$.\\

Berg, Saliola and Serrano \cite{BSS12} studied the \emph{Pieri operators} $D'_i$ acting on the affine nilCoxeter ring $\mathbb{A}_0$ defined by
\[ D'_i(A_w)= \sum_{w\longrightarrow_{(i)}u} A_u\]
where the sum is over all strong strips from $w$ to $u$ of length $i$. They proved that the operators $D'_i$ commute for all $i$ so that the symmetry of strong Schur functions follows. They also showed a series of properties of $D'_i$ that uniquely determine the operator $D'_i$.

\subsection{Cap operators}
Let us consider the affine type $A$ root datum and corresponding Kac-Moody flag variety $\widetilde{Fl}$ (See \cite{KK86} for details). Let $\{\xi^w | w\in W_\af\}$ denote the Schubert basis for the equivariant cohomology $H^*_T(\widetilde{Fl})$ of $\widetilde{Fl}$. One of main problems in affine Schubert calculus is to find a combinatorial formula for the structure constants $p_{u,v}^w$ where $ \xi^u \xi^v=\sum_w p_{u,v}^w\xi^w$. In \cite{KK86}, Kumar and Kostant identified the torus-equivariant cohomology $H^*_T(\widetilde{Fl})$ of the affine flag variety and the \emph{dual} of the nilHecke ring. Using this connection, one can compute the structure constants $p_{u,v}^w$ by calculating the coproduct structure constants of $\mathbb{A}$. More precisely, we have
\[\Delta(A_w)=\sum_{u,v\in W_\af} p_{u,v}^w A_u\otimes A_v. \]
For $u\in W_\af$, a \emph{cap operator} $D_u$ on $\mathbb{A}_0$ is defined by
$$D_u(A_w)=\sum_{\substack{v\in W_\af \\ \ell(v)=\ell(w)-\ell(u)}} p_{u,v}^w A_v.$$ 
Geometrically, the cap operator is the cap product on the ordinary homology and cohomology of the affine flag variety. More precisely, the cap operator $D_u$ can be considered as an element in $H^*(\widetilde{Fl})$ and $A_w$ can be considered in $H_*(\widetilde{Fl})$ so that the cap product $H^*(\widetilde{Fl})\times H_*(\widetilde{Fl})\rightarrow H_*(\widetilde{Fl})$ can be described by $(D_u,A_w)=D_u(A_w)$. Note that the cap operator only keep track of ordinary cohomology since $p_{u,v}^w$ is constant when $\ell(v)=\ell(w)-\ell(u)$. \\
 Let $\rho_i$ be the Pieri element $s_{i-1}\ldots s_1 s_0$ in $W_\af$ where indices are taken modulo $n$. We study the cap operators $D_{\rho_i}$ for $\rho_i$ and show that $D_{\rho_i}$ satisfy the properties of the Pieri operators $D'_i$ that Berg, Saliola and Serrano proved in \cite{BSS12}. This establishes the following main theorems equivalent to the affine Pieri rule conjectured in \cite{LLMS10}.

\begin{thm}\label{affinePieri}
For $w,u \in W_\af$ with $i=\ell(w)-\ell(u)\in \mathbb{N}$, $p^w_{\rho_i,u}$ counts the number of strong strips from $w$ to $u$.
\end{thm}
As a corollary, one can compute $p^w_{v,u}$ for $v=\rho_m^{-1}$ or $v=s_{i+m-1}s_{i+m-2}\cdots s_{i+1}s_i$ for any $i$ by applying automorphisms on Dynkin diagram of affine type $A$, namely, the map $j\mapsto -j$ and $j\mapsto j+i$ for $j\in \mathbb{Z}/n\mathbb{Z}$.\\
\smallskip

The paper is structured as follows. In section 2, we recall some notions about the affine symmetric groups, $k$-Schur functions, strong Schur functions
and root systems. In section 3, we define the affine nilHecke ring $\mathbb{A}$ and study its properties. In section 4, we recall some statements concerning the affine flag variety and its equivariant cohomology as well as the relationship between the equivariant cohomology of affine flag variety and the coproduct structure of the affine nilHecke ring. In section 5, we define the cap operator. We show that these cap operators for Pieri elements agree with the Pieri operators defined by Berg, Saliola and Serrano. We also prove the affine Pieri rule for the ordinary cohomology of the affine flag variety. In section 6, we apply the affine Pieri rule to show that the structure constants $p^w_{u,v}$ for the cohomology of the affine flag variety can be described in terms of strong Schur functions when $u$ is $0$-Grassmannian. This also gives a geometric interpretation of the skew strong Schur functions.
\section*{Acknowledgments}
I would like to thank Thomas Lam for helpful discussions and for introducing affine Schubert calculus to me. I also wish to thank KIAS for valuable support.
\section{Preliminaries}
\subsection{Affine symmetric group}

Positive integers $n\geq 2$ and $k=n-1$ will be fixed throughout the paper. Let $W_\af$ denote the affine symmetric group with simple generators $s_0,s_1,\ldots,s_{n-1}$ satisfying the relations
\begin{align*}
s_i^2&=1\\
s_is_{i+1}s_i&=s_{i+1}s_is_{i+1}\\
s_is_j&=s_js_i&& \text{if } i-j \neq 1,-1.
\end{align*}
where indices are taken modulo $n$. An element of the affine symmetric group may be written as a word in the generators $s_i$. A \emph{reduced word} of the element is a word of minimal length. The \emph{length} of $w$, denoted $\ell(w)$, is the number of generators in any reduced word of $w$. \\
The \emph{Bruhat order}, also called \emph{strong order}, on affine symmetric group elements is a partial order where $u<w$ if there is a reduced word for $u$ that is a subword of a reduced word for $w$. If $u<w$ and $\ell(u)=\ell(w)-1$, we write $u\lessdot w$. It is well-known that $u\lessdot w$ if and only if there exists a reflection $t\in \{ws_iw^{-1}|w\in W_\af, 0 \leq i<n \}$ such that $w=ut$ and $\ell(u)=\ell(w)-1$. For type $A$, the set $\{ws_iw^{-1}|w\in W_\af\}$ consists of transpositions $t_{ab}$. See \cite{BB05} for instance.\\
The subgroup of $W_\af$ generated by $\{s_1,\cdots,s_{n-1}\}$ is naturally isomorphic to the symmetric group $W$. The $0$-Grassmannian elements are minimal length coset representatives of $W_\af/ W$.

\subsection{Symmetric functions}
Let $\Lambda$ denote the ring of symmetric functions. For a partition $\lambda$, we let $m_\lambda,h_\lambda,e_\lambda,s_\lambda$ denote the monomial, homogeneous, elementary and Schur symmetric functions, respectively, indexed by $\lambda$. Each of these families forms a basis of $\Lambda$.\\
Let $\Lambda_{(k)}$ denote the subalgebra generated by $h_1,h_2,\ldots,h_k$. The elements $h_\lambda$ with $\lambda_1\leq k$ form a basis of $\Lambda_{(k)}$. Let $\Lambda^{(k)}=\Lambda/I_k$ denote the quotient of $\Lambda$ by the ideal $I_k$ generated by $m_\lambda$ with $\lambda_1>k$. The image of elements $m_\lambda$ with $\lambda_1\leq k$ form a basis of $\Lambda^{(k)}$. Lam \cite[Theorem 7.1]{Lam08} showed that $\Lambda_{(k)}$(resp. $\Lambda^{(k)}$) are isomorphic to homology(resp. cohomology) of affine Grassmannian as Hopf-algebras.\\

\subsection{$k$-Schur and Strong Schur functions}
There are many conjecturally equivalent definitions of $k$-Schur functions in \cite{LLMSSZ}. We follow the definition of $k$-Schur functions in \cite{LLMS10} using strong strips. Note that the marking of strong covers and strong strips defined in this extended abstract follows the notation of \cite{BSS11}. This differs from the definition of strong covers and strips in \cite{LLMS10} by reversing direction.\\

A \emph{marked strong cover} $C=\left(u\overset{a}\longrightarrow v \right)$ consists of $u,v\in W_\af$ and an integer $a$ such that $u=vt_{ij},v\lessdot u, a= u(j)=v(i)$ where $i\leq 0<j$. We use the notation inside$(C)=u$ and outside$(C)=v$.\\
A \emph{strong strip} S of length $i$ from $u$ to $v$, denoted by $u\longrightarrow_{(i)} v$, is a path
\[u\overset{a_1}\longrightarrow u_1 \overset{a_2}\longrightarrow\cdots \overset{a_i}\longrightarrow u_i=v\]
where $a_1>a_2>\ldots>a_i$. We let inside$(S)=u$, outside$(S)=v$, and size$(S)=i$.\\
A \emph{strong tableau} is a sequence $T=(S_1,S_2,\ldots)$ of strong strips $S_i$ such that outside$(S_j)=$ inside$(S_{j+1})$ for all $j \in \mathbb{Z}_{>0}$ and size$(S_i)=0$ for all sufficiently large $i$. We define inside$(T)=$ inside$(S_1)$ and outside$(T)=$ outside$(S_i)$ for $i$ large. The \emph{weight} $wt(T)$ of $T$ is the sequence 
\[wt(T)=(\text{size}(S_1),\text{size}(S_2),\ldots).\]
We say that $T$ has \emph{shape} $u/v$ where $u=$ inside$(T)$ and $v=$ outside$(T)$ so that we have $u>v$. 
\begin{definition}
For fixed  $u,v \in W_\af$, define the Strong Schur function
\[\emph{Strong}_{u/v}(x)=\sum_T x^{\text{wt}(T)}\]
where T runs over the strong tableaux of shape $u/v$.

\end{definition}

If $u$ is $0$-Grassmannian and $v$ is the identity element, $\textrm{Strong}_u(x)$ is a $k$-Schur function $s_{{\bf c}(u)}^{(k)}$ where ${\bf c}$ is a bijection between the set of 0-Grassmannian elements and the set of $k$-bounded partitions. For details, see \cite{LLMS10} for instance. \\

\subsection{Root systems}

We shall assume basic familiarity with Weyl groups, root systems, and weights. See \cite{Hum90} for details. \\
Let $A=(a_{ij})_{i,j\in I_\af}$ denote an affine Cartan matrix, where $I_{\af}=I\cup\{0\}$, so that $(a_{ij})_{i,j\in I}$ is corresponding finite Cartan matrix. For type $\widetilde{A}_{n-1}$(corresponding to $W_\af$), we have $I_{\text{af}}=\mathbb{Z}/n\mathbb{Z}$ and 
\[a_{ij}= \left\{\begin{array}{rl} 2 & \mbox{if } i=j   \\ -1 & \mbox{if }  |i-j|=1 \\ 0  &\mbox{otherwise}.\end{array} \right.\]
Let $R$ be the root system for $W$. Let $R^+, R^-$ denote the set of positive roots, negative roots respectively. Let $\{ \alpha_i\mid i\in I\}$ denote the simple roots and $\{\alpha_i^\vee\mid i\in I\}$ denote the simple coroots. Let $\theta$ denote the highest root $\alpha_1+\alpha_2+\ldots + \alpha_{n-1}$ of $R^+$. \\
Let $R_{\text{af}}$ and $R_{\af}^+$ denote the affine root system and positive affine roots. The positive simple affine roots (resp. coroots) are  $\{\alpha_i\mid i\in I_{\af}\}$ (resp. $\{\alpha_i^\vee\mid i\in I_{\af}\}$). The null root $\delta$ is given by $\delta=\alpha_0+\theta=\alpha_0+\cdots+ \alpha_{n-1}$. Similarily, the canonical central element $c$ is given by $\alpha_0^\vee+\cdots+ \alpha_{n-1}^\vee$. A root $\alpha$ is real if it is in $W_\af$-orbit of the simple affine roots, and imaginary otherwise. The imaginary roots are exactly $\{k\delta \mid k\in \mathbb{Z}\backslash \{0\}\}$. Every real affine root is of the form $\alpha+k\delta$, where $\alpha\in R$ and $k\in \mathbb{Z}$. The root $\alpha+k\delta$ is positive if $k>0$, or if $k=0$ and $\alpha\in R^+$. Let $R_{\text{re}}^+$ denote the set of positive roots in $R_\af$.\\
Let $Q=\oplus_{i\in I} \mathbb{Z}\cdot \alpha_i$ denote the root lattice and let $Q^\vee=\oplus _{i \in I} \mathbb{Z} \cdot \alpha_i^\vee$ denote the coroot lattice. Let $P$ and $P^\vee$ be the weight lattice and coweight lattice respectively. We have inclusions $Q \subset P, Q^\vee \subset P^\vee$ and a map $Q_{\af}=\oplus_{i \in I_\af} \mathbb{Z} \cdot \alpha_i \rightarrow P$ given by sending $\delta$ to $0$. Let $\langle\cdot,\cdot\rangle$ denote the pairing between $P$ and $P^\vee$ satisfying $\langle \alpha_i^\vee, \alpha_j\rangle=a_{ij}$.\\

Let $P_{\af}$ be the affine weight lattice and $\{\Lambda_i | i\in I_{\af} \}$ be the set of affine fundamental weights. Then 
\[P_{\af}= \mathbb{Z}\delta \oplus \bigoplus_{i\in I_{\af}} \mathbb{Z}\Lambda_i. \]

The affine Weyl group $W_\af$ acts on the affine weight lattice $P_\af$ and the affine coweight lattice by

\begin{align*}
s_i\cdot \lambda&=\lambda- \langle\alpha_i^\vee, \lambda\rangle \alpha_i\\
s_i \cdot \mu&=\mu - \langle\mu,\alpha_i\rangle\alpha_i^\vee.
\end{align*}
This is called the \emph{level zero action}. An element $w\in W_\af$ can be uniquely written of the form $ut_\mu$ where $u\in W$ and $\mu \in Q^\vee$. Then the level zero action on $P\subset P_\af$ is given by
$$ut_\mu\cdot \lambda = u \cdot \lambda.$$
For a real root $\alpha$, we let $s_\alpha$ denote the corresponding reflection, defined by $s_\alpha=w s_i w^{-1}$ if $\alpha=w\cdot\alpha_i$. The reflection $s_\alpha$ acts on weights by $s_{\alpha}\lambda=\lambda- \langle\alpha,\lambda\rangle\alpha$. For a strong cover $v\lessdot w$, let $\alpha_{v,w}$ denote the unique positive root satisfying the equation $v^{-1}w=s_{a_{v,w}}$.

\section{NilHecke algebra}
\subsection{Definition}
Let $S=\text{Sym}(P)$ be the polynomial ring having a variable for each free generator of $P$ and let $F=\text{Frac}(S)$ be the fraction field. \\

Define the $F$-vector space $F_{W_{\af}}=\bigoplus_{w\in W_{\af}} F w$ with basis $W_{\af}$, with product given by 
\[(pv)(qw)=\left( p(v\cdot q)\right) (vw) \quad \mbox{for } p,q \in F \mbox{ and } v,w \in W_{\af}.\]
For any real root $\alpha \in R_{\text{re}}$ define the element $A_\alpha\in F_{W_{\af}}$ by
\[A_\alpha= \alpha^{-1} (1-s_\alpha).\]
We write 
\[ A_i = A_{\alpha_i} \quad \mbox{for } i\in I.\]
For $\alpha=w \cdot \alpha_i \in R_{\text{re}}^+$, we have 
\[\begin{array}{rcl}
wA_i w^{-1} & = & A_\alpha  \\ s_\alpha A_\alpha & = &  A_\alpha \\ A_\alpha s_\alpha&=& - A_\alpha \\ A_\alpha^2&=&0. \end{array}\]
The $A_i$ satisfy the braid relations as the $s_i$ in $W_{\af}$, i.e., $A_iA_{i+1}A_i=A_{i+1}A_iA_{i+1}$. Therefore it makes sense to define
\[ \begin{array} {rlll} A_w&=&A_{i_1}\cdots A_{i_l} & \mbox{where} \\ w&=&s_{i_1}\cdots s_{i_l} & \mbox{is a reduced decomposition.} \end{array} \]
One can check that 
\[A_vA_w=\left\{ \begin{array}{ll} A_{vw} & \mbox{if } \ell(vw)=\ell(v)+\ell(w) \\ 0 & \mbox{otherwise.}\end{array} \right. \]
The \emph{nilCoxeter algebra} $\mathbb{A}_0$ is the subring of $F_{W_{\af}}$ generated by $A_i$ over $\mathbb{Z}$ for $i\in I_{\af}$. The set $\{A_w|w\in W_\af\}$ forms a basis of $\mathbb{A}_0$ over $\mathbb{Z}$. The \emph{nilHecke algebra} $\mathbb{A}$ is the subring of $F_{W_{\af}}$ generated by $S$ and $\mathbb{A}_0$.

We state the following result \cite[Proposition 4.30]{KK86}. Let $w\in W_\af$ and $\lambda\in P$. Then
\begin{align} A_w\lambda= (w\cdot \lambda) A_w + \sum_{ws_\alpha:\ell(ws_\alpha)=\ell(w)-1} \langle\alpha^\vee, \lambda\rangle A_{ws_\alpha} \end{align}
where $\alpha$ is always taken to be a positive root in $R_\af^+$. The coefficients $\langle\alpha^\vee,\lambda\rangle$ are known as \emph{Chevalley coefficients}.

\begin{rmk}
The affine nilHecke algebra $\mathbb{A}$ defined above is slightly different from those defined in \cite{KK86}. In \cite{KK86}, Kostant and Kumar set $S=\text{Sym}(P_{\text{af}})$ instead of $\text{Sym}(P)$. We have a projection $P_{\text{af}} \rightarrow P$ and the map induces $\text{Sym}(P_{\text{af}})\rightarrow \text{Sym}(P)$ which is compatible with Theorem \ref{cohomology}. The image of the null root $\delta=\alpha_0+\theta$ in $P$ is 0. Geometric interpretation of $P$ and $P_\af$ will be discussed in section \ref{affineflagvariety}.
\end{rmk}

\subsection{Cyclically decreasing elements and Peterson algebra}
A word $s_{i_1}s_{i_2}\cdots s_{i_l}$ with indices in $\mathbb{Z}/n\mathbb{Z}$ is called \emph{cyclically decreasing} if each letter occurs at most once and whenever $s_i$ and $s_{i+1}$ both occurs in the word, $s_{i+1}$ precedes $s_i$. For $J\varsubsetneq\mathbb{Z}/n\mathbb{Z}$ there is the unique cyclically decreasing element $w_J$ with letters $\{s_j|j\in J\}$. For $i \in \{ 0,1,\ldots,n-1\}$, let
\[\h_i=\sum\limits_{\substack{ J\subset I_{\af} \\ |J|=i }} A_{w_J} \in \mathbb{A}_0 \]
where $\h_0=1$ and $\h_i=0$ for $i< 0$ by convention. Lam \cite{Lam06} showed that the elements $\{\h_i\}_{i<n}$ commute and freely generate a subalgebra $\mathbb{B}$ of $\mathbb{A}_0$ called \emph{affine Fomin-Stanley algebra}. It is well-known that $\mathbb{B}$ is isomorphic to $\Lambda_{(k)}$ via the map sending $\h_i$ to $h_i$. Therefore, the set $\{\h_\lambda=\h_{\lambda_1}\ldots \h_{\lambda_l} | \lambda_1\leq k\}$ forms a basis of $\mathbb{B}$.
\\

\subsection{Coproduct structure of $\mathbb{A}$} \label{coproduct}
There is a coproduct structure on the nilHecke ring $\mathbb{A}$. Kostant and Kumar showed that the structure constants of the coproduct in $\mathbb{A}$ are the same as the structure constants of the equivariant cohomology ring of $\widetilde{Fl}$. Details and proofs can be found in \cite{KK86,LLMSSZ,Pet97}.\\
Let $\Delta:F_{W_\af}\rightarrow F_{W_\af}\otimes_F F_{W_\af}$ be the left $F$-linear map defined by
\[\Delta(w)=w \otimes w \quad \mbox{for all } w\in W_\af\]

\begin{thm} \label{delta}
The map $\Delta:F_{W_\af}\rightarrow F_{W_\af}\otimes_F F_{W_\af}$ induces the unique left $S$-module homomorphism $\Delta: \mathbb{A} \rightarrow \mathbb{A}\otimes_S \mathbb{A}$ such that

\[\begin{array}{rll}
\Delta(A_i)&=A_i\otimes 1 + s_i \otimes A_i & \\
 &= 1\otimes A_i + A_i \otimes s_i & \text{ for all } i \in I \\
\Delta(ab)&=\Delta(a)\Delta(b) & \text{ for all } a,b\in \mathbb{A}. \\
\end{array}
\]
\end{thm}
For $a,b \in F_{W_\af}$ and any expression of $\Delta(a)$ and $\Delta(b)$ of the form
\begin{align*} \Delta(a)&=\sum_{a} a_{(1)}\otimes a_{(2)}\\ \Delta(b)&=\sum_{b} b_{(1)}\otimes b_{(2)}. \end{align*} 
The product of $\Delta(a)$ and $\Delta(b)$ can be computed by the naive componentwise product
\begin{align}\label{naive}
\Delta(ab)=\Delta(a)\Delta(b)= \sum_{(a),(b)} a_{(1)}b_{(1)}\otimes a_{(2)}b_{(2)}.
\end{align}

For $w,u_1,u_2\in W_\af$, the \emph{equivariant Schubert structure constants} $p^w_{u_1,u_2}\in S$ are defined as the coefficients in the expansion $\Delta(A_w)=\sum p^w_{u_1,u_2} A_{u_1} \otimes A_{u_2}$. The followings are properties of $p^w_{u_1,u_2}$.

\begin{thm} \ \begin{enumerate}
\item $p^w_{u_1,u_2}=0$ unless $w\geq u_1$ and $w\geq u_2$.\\
\item $p^w_{u_1,u_2}$ is homogeneous of degree $\ell(u_1)+\ell(u_2)-\ell(w)$. \\
\item \cite{Gra01,Kum02} $(-1)^{\ell(u_1)+\ell(u_2)-\ell(w)}p^w_{u_1,u_2}\in \mathbb{Z}_{\geq 0} [\alpha_i | i\in I_\af].$
\end{enumerate} \label{essc}
\end{thm}
The last property in Theorem \ref{essc} is called \emph{Graham positivity}. 
One can use Theorem \ref{delta} to compute $p^w_{u_1,u_2}$ explicitly.
\begin{thm}\label{p}
\cite{KK86}
Let $w=s_{i_1}s_{i_2}\ldots s_{i_l}$ and $l=\ell(w)$. Then
$$p^w_{u_1,u_2} = \sum\limits_{\substack{ A_{s_{i_{j_1}}} A_{s_{i_{j_2}}}\ldots A_{s_{i_{j_k}}} = A_{u_1} \\
 j_1<j_2<\ldots<j_k }}\prod_{p=1}^l f_{j_1,j_2,\ldots,j_k}(A_{s_{i_{p}}})\Bigg|_{A_{u_2}}$$
where $f_{j_1,j_2,\ldots,j_k}(A_{s_{i_{p}}})= s_{i_{p}}$ if $p \in \{j_1,j_2,\ldots,j_k\}$ and $=A_{s_{i_{p}}}$ otherwise, and $\sum c_w A_w |_{A_u}=c_u$.

\end{thm}

\section{Affine flag varieties} \label{affineflagvariety}
In this section, we define the Kac-Moody flag variety $\widetilde{Fl}$ and establish the relationship between the equivariant cohomology of $\widetilde{Fl}$ and the coproduct structure on $\mathbb{A}$. There are two definitions of Kac-Moody flag variety $\widetilde{Fl}$ in \cite{LLMSSZ}, but we only recall $\widetilde{Fl}$ as the Kac-Moody flag ind-variety in \cite{KK86,Kum02}. \\
Let $G_\af$ denote the Kac-Moody group of affine type associated with $SL(n)$ and let $B_\af$ denote its Borel subgroup. The Kac-Moody flag ind-variety $\fl=G_\af/B_\af$ is paved by cells $B_\af\dot{w}B_\af/B_\af\cong\mathbb{C}^{\ell(w)}$ whose closure $X_w$ is called the Schubert variety. Schubert variety defines a Schubert class $[X^w]_{T_\af} \in H^*_{T_\af}(\widetilde{Fl})$ where $T_\af$ is the maximal torus in $B_\af$. For our setting, we consider $H^*_{T}(\widetilde{Fl})$ where $T$ is a maximal torus in $SL(n)$. Kumar and Kostant \cite{KK86} identified the equivariant cohomology ring of $\widetilde{Fl}$ with the \emph{dual} of nilHecke ring $\mathbb{A}$. By restricting the group action from $T_\af$ to $T$, this identification shows that the structure constant of $H^*_T(\widetilde{Fl})$ is the same as $p_{u,v}^w$ defined in Section \ref{coproduct}.

\begin{thm} \label{cohomology}
\cite{KK86} The $T$-equivariant cohomology of $\widetilde{Fl}$ has a basis $\{[X^w]_T \in H^*_T(\widetilde{Fl})\}$ over $S\cong H^*_T(pt)$. Moreover, the structure constants of $H^*_T(\widetilde{Fl})$ are $p^w_{u_1,u_2}$, i.e.,
$$ [X^{u_1}]_T [X^{u_2}]_T=\sum_w p^w_{u_1,u_2}[X^w]_T. $$
\end{thm}
\begin{rmk}\label{ordinary}

The ordinary cohomology of $\widetilde{Fl}$ has a basis $\{[X^w] \in H^*(\widetilde{Fl})\}$ over $\mathbb{Z}$ and corresponding structure constants are just $\phi(p^w_{u_1,u_2})$ where $\phi:S\rightarrow\mathbb{Z}$ is the evaluation at $0$.
\end{rmk}

\section{Pieri operators and Affine Pieri rule for type $A$}

Let $\phi: S\rightarrow \mathbb{Z}$ denote the map sending a polynomial to its constant term, called the \emph{evaluation map} at $0$. It extends to the map $\phi: \mathbb{A}\rightarrow\mathbb{A}_0$ given by $\phi(\sum_w a_wA_w)=\phi(a_w)A_w$. \emph{Pieri elements} $\rho_i$ are defined by $s_{i-1}s_{i-2}\ldots s_{1}s_0$ for $i>0$.
For $w\in W_{\text{af}}$, let us define \emph{the cap operator} $D_{w}$ on $\mathbb{A}$ (and $\mathbb{A}_0$) by

\[D_w(A_v):= \phi\left(\sum_{u\in W_\af} p^v_{w,u} A_u\right) =\sum_{\substack{u\in W_\af\\\ell(u)=\ell(v)-\ell(w) }} p^v_{w,u} A_u.  \]
 If $w=\rho_i$, we call $D_{\rho_i}$ the \emph{Pieri operator} and denote by $D_i$.

Note that $p^v_{w,u}\in \mathbb{Z}_{\geq 0}$ when $\ell(v)=\ell(w)+\ell(u)$ by Theorem \ref{essc}. Finding a combinatorial formula for such $p^v_{w,u}$ is one of the important problems in Schubert calculus and it is not completely known even for finite flag varieties. In this section, we will prove that $D_{i}$ is the same as the Pieri operator $D'_i$ defined by Berg, Saliola, Serrano in \cite{BSS12}. This identification gives a combinatorial description of the Pieri rule for the ordinary cohomology of $\widetilde{Fl}$ and therefore proves the main Theorem \ref{affinePieri}.\\
Let us state properties of $D'_i$ described in \cite{BSS12}.

\begin{thm} \cite[Theorem 4.8]{BSS12}\label{bss1}
Suppose $w \in W_{\text{af}}$ and $v \in W$. Then
\[D'_i (A_w A_v)=D'_i (A_w) A_v. \]

\end{thm}

\begin{lemma}\label{bss2} \cite[Lemma 4.5]{BSS12} For $r<n$ and $i\geq 1$,
\[D'_i(\h_r)=\h_{r-i}.\]
\end{lemma}

\begin{thm} \cite[Proposition 4.3]{BSS12} \label{bss3} For $p < n$,$w\in W_\af$ and $i\geq 1$, we have
\[ D'_i(\h_p A_w) = \sum_{j= 0}^i D'_j(\h_{p}) D'_{i-j}(A_w) =\sum_{j= 0}^i \h_{p-j} D'_{i-j}(A_w).\]
\end{thm}

For $w\in W_\af$, there is a unique decompostion $w=w_\lambda w(0)$  where $w_\lambda$ is $0$-Grassmannian and $w(0)\in W$. Note that $\lambda={\bf c}(w_\lambda)$ where ${\bf c}$ is the bijection between $0$-Grassmannian elements and $k$-bounded partitions mentioned in section 2. In this setup, the set $\{ \h_\lambda A_{w(0)} : w= w_\lambda w(0) \in W_{\text{af}} \}$ forms a basis of $\mathbb{A}$ (See \cite{BSS12}). Therefore, above theorems uniquely determine $D'_i$. In this section, we prove the same theorems for $D_{i}$. This implies that $D_i=D'_i$ and Theorem \ref{affinePieri} follows. 

\begin{thm}\label{0}
Suppose $w \in W_{\text{af}}$ and $v \in W$. Then
\[D_{i} (A_w A_v)=D_{i}(A_w) A_v \]

\end{thm}

{\it Proof.} Let $F(A_w;j_1,j_2,\ldots,j_k) = \prod_{p=1}^l f_{j_1,j_2,\ldots,j_k}(A_{s_{i_{p}}})$ for $w=s_{i_1}s_{i_2}\ldots s_{i_l}$ and $1 \leq j_1<j_2<\ldots<j_k\leq l$.
It is enough to show that
$$ F(A_wA_v;j_1,j_2,\ldots,j_k) =F(A_w;j_1,j_2,\ldots,j_k) A_v $$
for all $1\leq j_1<j_2<\ldots <j_k \leq l$ by Theorem \ref{p}. Note that any reducd word of $v$ does not contain $s_0$ since $v\in W$, and $\rho_i$ has the unique reduced word ending with $s_0$, namely $s_i s_{i-1} \ldots s_1 s_0$. Let $w=s_{i_1}s_{i_2}\ldots s_{i_{\ell(w)}}$ and $v=s_{a_1}s_{a_2}\ldots s_{a_{\ell(v)}}$. Then $a_b\neq 0$ for $1\leq b\leq \ell(v)$. Therefore, we have

\begin{align*}
 F(A_wA_v;j_1,j_2,\ldots,j_k) &= \prod_{p=1}^{\ell(w)} f_{j_1,j_2,\ldots,j_k}(A_{s_{i_{p}}} A_v) \\
&= \prod_{p=1}^{\ell(w)} f_{j_1,j_2,\ldots,j_k}(A_{s_{i_{p}}}) A_v \\
&=F(A_w;j_1,j_2,\ldots,j_k) A_v 
\end{align*}
and the theorem follows. \qed\\

\begin{thm}\label{dh}
For $r<n$ and all $i\geq 1$,
\[D_i(\h_r)=\h_{r-i}.\]
\end{thm}

{\it Proof.} This theorem follows as a corollary of \cite[Lemma 7.7]{Lam08}.
\begin{lemma}\cite[Lemma 7.7]{Lam08} For $r<n$,
$$\phi(\Delta(\h_r))=\sum_{0\leq j\leq r} \h_{r-j} \otimes \h_{j}.$$
\end{lemma}
Since $\rho_i$ appears in $\h_i$, terms appearing in $\phi(\Delta(\h_r))$ with second entry $\rho_i$ are exactly $\h_{r-i}\otimes \rho_i$.
\qed \\

In order to prove analogue of Theorem \ref{bss3} for $D_i$, we need the following lemmas. 

\begin{lem}\label{hp} \cite[Lemma 7.2]{Lam08}
Let $b\in \mathbb{B}$ and $s\in S$. Then
$$\phi(b \cdot s )= \phi(s)b=b \cdot\phi(s).$$

\end{lem}

\begin{lem} \label{translation}
Let $\psi: W_\af \rightarrow W_\af$ be the group automorphism defined by $s_i\mapsto s_{i-1}$. The automorphism induces the ring automorphism $\psi: \mathbb{A}_0 \rightarrow \mathbb{A}_0$. For $b\in \mathbb{B}$, we have
$$\psi(b)=b.$$
\end{lem}
{\it Proof.} Since $\mathbb{B}$ is generated by $\h_i$ for $i<n$ and the equality $\phi(\h_i)=\h_i$ is obvious from the definition of $\h_i$, we are done.
\qed \\
 
Now we are ready to prove the following property of $D_i$.
\begin{thm}\label{dual} For $p < n$,$w\in W_\af$ and $i\geq 1$, we have
\begin{displaymath}D_{i}(\h_p A_w)= \sum_{j=0}^i D_{j} (\h_p) D_{i-j}(A_w).\end{displaymath}
\end{thm}

{\it Proof.} By the definition of the Pieri operator and Remark \ref{ordinary}, we have
\begin{align*} 
D_{\rho_i}(\h_p A_w) &= \sum_{|S|=p} D_{\rho_i}(A_{w_s}A_w) \\
&=\sum_{|S|=p}\sum\limits_{\substack{u}}\phi\left(\text{coefficient of } A_u\otimes A_{\rho_i} \text{ in } \Delta(A_{w_s}A_w)\right)A_u \\
&=\sum_{|S|=p}\sum_{\substack{u}}\phi\left(\text{coefficient of } A_u\otimes A_{\rho_i} \text{ in } \Delta(A_{w_s})\Delta(A_w)\right)A_u.\\ 
\end{align*}
Note that the last equality follows from Theorem \ref{delta}. By (\ref{naive}), we can compute $ \Delta(A_{w_s})$ and $\Delta(A_w)$ seperately and multiply them together. Let $\widetilde{\rho}_{ij}$ denote the element $\rho_i \rho_{i-j}^{-1}=s_{i-1}s_{i-2}\ldots s_{i-j}$ for $i\geq j$. Then $\psi^{i-j}(\widetilde{\rho}_{ij})=\rho_j$. Therefore, we have
\begin{align*}
D_{\rho_i}(\h_p A_w)&=\phi\Big(\sum_{|S|=p}\sum_{j=0}^i\sum_{\substack{u_1,u_2 }} \left(\text{coefficient of } A_{u_1}\otimes A_{\widetilde{\rho}_{ij}} \text{ in }\Delta(A_{w_s})\right)A_{u_1}  \\
&\cdot \left(\text{coefficient of } A_{u_2}\otimes A_{\rho_{i-j}} \text{ in } \Delta(A_{w})\right) A_{u_2}\Big) \\
&=\phi\Big(\sum_{j=0}^i\sum_{\substack{u_1,u_2 }}\sum_{|S|=p} p_{u_1,\widetilde{\rho}_{ij}}^{w_s}A_{u_1}\cdot p_{u_2,\rho_{i-j}}^{w} A_{u_2}\Big) \\
&=\phi\Big(\sum_{j=0}^i\sum_{\substack{u_1,u_2 }}\sum_{|S|=p} \phi\left(p_{u_1,\widetilde{\rho}_{ij}}^{w_s}\right)A_{u_1}\cdot p_{u_2,\rho_j}^{w} A_{u_2}\Big) \\
&=\sum_{j=0}^i\sum_{\substack{u_2 }}\phi\Big(D_{\widetilde{\rho}_{ij}}(\h_{p})\cdot p_{u_2,\rho_j}^{w}\Big) A_{u_2} \\
&=\sum_{j=0}^i\sum_{\substack{u_2 }}\phi\Big(D_j(\h_{p})\cdot p_{u_2,\rho_j}^{w}\Big) A_{u_2} \mbox{ (By Lemma \ref{translation})}\\
&=\sum_{j=0}^i\sum_{\substack{u_2 }}\phi\Big(\h_{p-j}\cdot p_{u_2,\rho_j}^{w}\Big) A_{u_2} \mbox{ (By Theorem \ref{dh})}\\
&=\sum_{j=0}^i\sum_{\substack{u_2}}\h_{p-j}\cdot \phi\left(p_{u_2,\rho_j}^{w}\right) A_{u_2} \mbox{ (By Lemma \ref{hp} with $b=\h_{p-j}$)}\\
&=\sum_{j=0}^i \h_{p-j} D_{\rho_{i-j}}(A_w).\qed
\end{align*}
\\
{\it Proof of Theorem \ref{affinePieri}.} Theorem \ref{0}, \ref{dh}, \ref{dual} uniquely determine the Pieri operator $D_i$ and we have $D_i=D'_i$ by Theorem \ref{bss1}, \ref{bss2}, \ref{bss3}. Therefore, Theorem \ref{affinePieri} follows.\qed\\
\begin{rmk}
Theorem \ref{dual} generalizes \cite[Theorem 3.8]{BSS11}. In fact, the author believes that there is a gap in the first few lines in their proof of \cite[Theorem 3.8]{BSS11}, and Theorem \ref{dual} resolves the problem as well.
\end{rmk}

\section{Applications}
Let $W^0$ denote the set of $0$-Grassmannian elements in $W_\af$. In this section, we prove that Theorem \ref{affinePieri} implies the following conjecture, as stated in \cite{LLMS10}.

\begin{thm}\cite[Conjecture 4.18]{LLMS10}\label{strong}
Let $w,u\in W_\af$ be two affine permutations. Then
\begin{enumerate}
\item $\Strong_{w/u}(x)\in \Lambda$.
\item $\Strong_{w/u}(x)\in \Lambda_{(n)}$.
\item $\Strong_{w/u}(x)=\sum_{v \in W^0} p^w_{u,v} \Strong_v(x)$.
\end{enumerate}
\end{thm}
Note that Berg, Saliola and Serrano \cite{BSS12} proved $(1),(2)$ of Theorem \ref{strong} and $(3)$ when $w,v$ are in $W^0$. To be more precise, they showed that for $w,v,u \in W^0$, $p^w_{u,v}$ is the coefficient of $\tilde F_w$ in $\tilde F_u \tilde F_v$, where $\tilde F_w$ is the \emph{affine Stanley symmetric function} labeled by $w$. Lam, Lapointe, Morse and Shimozono \cite{LLMS10} introduced \emph{weak Schur functions} defined by
$$ \Weak_{w/u}(x)= \sum_\lambda \langle \h_\lambda A_u , A_w \rangle_{\mathbb{A}_0} m_\lambda(x)$$
for $w,u \in W_\af$ where the summation is over all $k$-bounded partitions and $\langle \cdot, \cdot \rangle_{\mathbb{A}_0}$ is the inner product on $\mathbb{A}_0$ satisfying $\langle A_{w_1}, A_{w_2} \rangle_{\mathbb{A}_0} = \delta_{w_1, w_2}$. If $w \in W^0$ and $v$ is the identity element, $\Weak_{w/u}(x)$ is the affine Stanley symmetric function $\tilde F_w$. \\
Let $Gr$ denote the affine Grassmannian associated with $SL(n,\mathbb{C})$. Since $Gr \cong G_\af/\mathcal{P}$ for the affine Kac-Moody group $G_\af$ and its maximal parabolic subgroup $\mathcal{P}$, we have the Schubert bases
$$\{\zeta^w\in H^*(Gr,\mathbb{Z}) \mid w \in W^0 \}$$
$$\{\zeta_w\in H_*(Gr,\mathbb{Z}) \mid w \in W^0 \}.$$
Lam \cite{Lam08} showed that $\zeta_w$ and $\zeta^w$ can be represented by $k$-Schur functions $\Strong_w(x)$ and affine Stanley symmetric functions $\Weak_w(x)$ via the isomorphisms $H_*(Gr)\cong \Lambda_{(n)}$ and $H^*(Gr)\cong \Lambda^{(n)}$. In fact, \emph{affine insertion} \cite{LLMS10} offers a duality between weak and strong orders combinatorially.\\
Define the \emph{affine Cauchy kernel} $\Omega_n(x,y)$ by
\begin{align*}
\Omega_n(x,y)&= \prod_i \left( 1+ y_i h_1(x) + y_i^2 h_2(x) + \cdots + y_i^{n-1} h_{n-1}(x)\right)\\
&=\sum_{\lambda : \lambda_1<n} h_\lambda(x) m_\lambda(y)\\
&=\sum_{\lambda} m_\lambda(x) \tilde{h}_\lambda(y)
\end{align*}
as an element in $\Lambda_{(n)}(x)\widehat{\otimes}\Lambda^{(n)}(y)\subset\Lambda(x)\widehat{\otimes}\Lambda^{(n)}(y)$ where $\tilde{h}_\lambda(y)$ is the image of $h_\lambda(y)$ in $\Lambda^{(n)}$. The second equality follows from the well-known equality $\sum_\lambda h_\lambda(x)m_\lambda(y)=\sum_\lambda m_\lambda(x)h_\lambda(y)$ after taking the projection map $\Lambda(x)\widehat\otimes\Lambda(y)\rightarrow\Lambda(x)\widehat\otimes\Lambda^{(n)}(y)$.\\
The duality between strong and weak orders produces the following affine Cauchy identity \cite[Corollary 4.6]{LLMS10}.

\begin{cor} (Affine Cauchy Identity) \label{affinecauchy}
The following identity holds in the formal power series ring $\mathbb{Z}[[x_1,x_2,\ldots,y_1,y_2,\ldots]]$:
$$\Omega_n(x,y)=\sum_{w\in W^0} \Strong_w(x)\Weak_w(y).$$
\end{cor}
Now we are ready to prove Theorem \ref{strong}, following \cite[Proposition 4.19]{LLMS10}. \\
{\it Proof of Theorem \ref{strong}.}
Let $\{\xi^w\in H^*(\widetilde{Fl}) \mid w \in W \}$ $(\{\xi_w\in H_*(\widetilde{Fl}) \mid w \in W \})$ denote the Schubert bases for the cohomology (resp. homology) of the affine flag variety. In this setting, Theorem \ref{affinePieri} can be described as follows:
\begin{thm}\label{affinePieri2}
Let $w\in W_\af$ and $1\leq m$. Then in $H^*(\widetilde{Fl})$ we have
$$ \xi^{\rho_m} \xi^w = \sum_S \xi^{\outside(S)},$$
where the sum runs over strong strips $S$ of size $m$ such that $\inside(S)=w$.
\end{thm}
Theorem \ref{affinePieri2} completely determines the action of $H^*(Gr)$ on $H^*(\widetilde{Fl})$, obtained from the inclusion $H^*(Gr) \subset H^*(\widetilde{Fl})$. \\
We may consider the affine Cauchy kernel $\Omega_n$ as an element of the completion $\Lambda \widehat{\otimes} H^*(Gr)$ via the isomorphism $\Lambda^{(n)}\cong H^*(Gr)$. Let $\langle \cdot, \cdot \rangle_{\widetilde{Fl}}$ denote the inner product on $H^*(\widetilde{Fl})$ defined by $\langle \xi^w, \xi^u \rangle=\delta_{wu}$. By the definition of $\Omega_n$ and Theorem \ref{affinePieri2}, we have
$$ \Strong_{w/u}(x)= \langle \Omega_n \cdot \xi^u, \xi^w \rangle_{\widetilde{Fl}}$$
where $\Omega_n\cdot \xi^u \in \Lambda\widehat{\otimes} H^*(Gr)$. Note that we have also used the definition of the strong Schur functions and the equality $\Weak_{\rho_i}(y)=\tilde{h}_i(y)$ (See \cite{LLMS10}). By Corollary \ref{affinecauchy}, we may also write
$$ \Omega_n=\sum_{v\in W^0} \Strong_v(x)\otimes \xi^v$$
so that 
\begin{align*}
 \Strong_{w/u}(x)&= \sum_{v \in W^0} \Strong_v(x) \langle \xi^v\xi^u,\xi^w\rangle_{\widetilde{Fl}}\\
&=\sum_{v \in W^0} p^w_{u,v} \Strong_v(x).
\end{align*}
Since $\Strong_v(x)$ are $k$-Schur functions in $\Lambda_{(n)}$, we have $\Strong_{w/u}\in \Lambda_{(n)}$. \qed\\
\begin{rmk}
The skew strong Schur function $\Strong_{w/u}(x)$ can be geometrically interpreted as the image of $D_u(\xi_w)$ under the projection $H_*(\widetilde{Fl})\rightarrow H_*(Gr)$ where $D_u$ is the cap product acting on $H_*(\widetilde{Fl})$. Indeed, by the definition of $D_u$ and Remark \ref{ordinary} we have $D_u(\xi_w)=\sum_v p^w_{uv}\xi_v$. Since the projection $H_*(\widetilde{Fl})\rightarrow H_*(Gr)$ maps $\xi_w$ to $\Strong_w(x)$, it sends $D_u(\xi_w)$ to $Strong_{w/u}(x)$ by Theorem \ref{strong}(3). Therefore, $\Strong_{w/u}(x)$ captures the $0$-Grassmannian part of $D_u(\xi_w)$. Note that when $w$ is not $0$-Grassmannian, we have $\Strong_w(x)=0$ (See \cite{LLMS10}).
\end{rmk}
Note that Lam, Lapointe, Morse and Shimozono \cite{LLMS10} defined the marked strong cover $wt_{ij} \rightarrow w$ where $i \leq l < j$ for a fixed constant $l$ instead of $i \leq 0 <j$ in our setup. The affine Pieri rule conjectured in \cite{LLMS10} is the following.
\begin{thm}\cite[Conjecture 4.15]{LLMS10}\label{affinePieri3}
Let $w\in W_\af$ and $1\leq m$. Then in $H^*(\widetilde{Fl})$ we have
$$ \xi^{\psi^{-l}(\rho_m)} \xi^w = \sum_S \xi^{\outside(S)},$$
where the sum runs over strong strips $S$ of size $m$ with respect to $l$ such that $\inside(S)=w$.
\end{thm}
We can obtain Theorem \ref{affinePieri3} by applying $\psi^l$ to Theorem \ref{affinePieri2}. Recall that $\omega$ is the automorphism of $W_\af$ sending $s_i$ to $s_{n-i}$. By applying $\omega$ to Theorem \ref{affinePieri2}, we get the dual affine Pieri rule.
\begin{thm}
Let $w\in W_\af$ and $1\leq m$. Then in $H^*(\widetilde{Fl})$ we have
$$ \xi^{\omega(\rho_m)} \xi^w = \sum_S \xi^{\omega(\outside(S))},$$
where the sum runs over strong strips $S$ of size $m$ such that $\inside(S)=\omega(w)$.
\end{thm}

Since the image of $\xi^{\omega(\rho_m)}$ under the projection $H^*(\widetilde{Fl})\rightarrow H^*(Gr)\cong \Lambda^{(n)}$ is  $\Weak_{\omega(\rho_m)}=\tilde{e}_m$, the image of $e_m\in \Lambda^{(n)}$, so we recover the dual Pieri rule for the affine Grassmannian in \cite{LLMS10}.

\end{document}